\documentclass[10pt]{article}
\usepackage{amsmath}
\usepackage{amsfonts}
\usepackage{amssymb}
\usepackage{amsthm}
\newtheorem{thm}{Theorem}

\newtheorem{lem}{Lemma}

\theoremstyle{definition}
\newtheorem{defn}{Definition}
\newcommand{\Cal}{{\rm Cal}}
\newcommand{\Sol}{{\rm Sol}}
\newcommand{\Met}{{\rm Met}}
\newcommand{\Z}{{\bf Z}}
\newcommand{\R}{{\bf R}}
\newcommand{\E}{{\bf E}}
\newcommand{\p}{{\bf p}}
\author{A.~M.~Vershik\thanks{%
Supported in part by the RFBR grant
02-01-00093 and by the Russian President grant
for Support of Leading Scientific Schools NSh.2251.2003.1}
\and S.~V.~Dobrynin\footnotemark[1]}
\title{Geometrical approach to the free solvable groups}
\date{20.09.04}
\begin{document}

\maketitle

{\abstract We give a topological interpretation of the free metabelian group, following the plan
described in \cite{Ver1},\cite{Ver2}. Namely we represent the free metabelian group with
$d$-generators as extension the group of the first homology of the $d$-dimensional lattice as
Cayley graph of the group ${\Bbb Z}^d$ with a canonical 2-cocycle. This construction open the
possibility to study metabelian groups from new points of view; in particular to give a useful
normal forms of the elements of the group, applications to the random walks and so on. We also
describe the satellite groups which correspond to all 2-cocycles of cohomology group associated
with the free solvable groups. The homology of the Cayley graph can be used for the studying of the
wide class of groups which including the class of all solvable groups.}

\section{Introduction}
In this paper we suggest a useful topological model of free solvable  groups of given degree of
solvability (mainly of degree two, i.e., the  metabelian free group). Using simple topological
tools, we obtain a recursive (with respect to the degree of solvability) description of free
solvable groups. The idea (see \cite{Ver1, Ver2}) is to consider the Cayley graph of a finitely
generated group {\it as a one-dimensional complex} and use classical topological results. We start
from the general situation. The Cayley graph $Cal(G)$ of any countable group $G$ as one-dimensional
complex is homotopic to a bouquet of 1-spheres, and hence its fundamental group either is trivial
(for the Cayley graph of a free group with standard generators), or is a free group. On the other
hand, if a given group $G$ is finitely generated and represented as a quotient of a free group
$F_d$, the corresponding normal subgroup $N$ of the free group is also free, i.e., isomorphic to
the fundamental group (with identity as a basepoint) of the Cayley graph of the group $G$. The free
group $F_d$ acts on the fundamental group of Cayley graph $\pi_1(Cal(G))$ as well as on the normal
subgroup $N$ (by inner automorphisms). The first step is almost tautological: there is the natural
isomorphisms between fundamental group $\pi_1(Cal(G))$ of the Cayley graph of the group and the
normal subgroup $N$ which is equivariant with respect to those actions of free group $F_d$ (Theorem
1 -paragraph 2).

The next step is more important and concerned to the group of the first homology of Cayley group.
The classical Gurevich's theorem asserts that the first homology group is the quotient of the
fundamental group over its commutant $H_1(Cal(G))=(\pi_1(Cal(G))/(\pi_1(Cal(G))'$. From other side
the Abelianizing the normal subgroup $N$ (i.e., factorizing over its commutant $N/N'$) allows the
to define the action of the group itself $G=(F_d)/N$ on the group $N/N'$ and in the same time the
group $G$ acts on the its graph Cayley $Cal(G)$ and consequently on the first homology group of the
Cayley graph $H_1(Cal(G))$; the theorem 2 asserts that there is natural  $G$-equivariant
isomorphism between $N/N'$ and $H_1(Cal(G))$ (Theorem 2).

The third step gives us the construction of the extensions. The interpretation of the quotient
$(F_d)/N'$ of the free group $F_d$ over the commutant $N'$ of the normal subgroup $N$ as an
extension of the homology group $H_1(Cal(G))$ of the Cayley graph by the group $G=(F_d)/N$. For
this we must define the canonical 2-cocycle of the group $G$ with values in the the group
$H_1(Cal(G))$. Suppose, the normal subgroup $N=(F_d)^n$ is the $n$-th commutant of the free group
$F_d$, then the quotient group is the free solvable group of degree $n$ with $d$ generators
$Sol(n,d)$; In order to give the recursive interpretation of it we define a canonical 2-cocycle of
group of previous degree of solvability $Sol(n-1,d)$ with values in the group of homology
$H_1(Cal(Sol(n-1,d)))$ of its Cayley graph. We obtain (Theorem 3) a recursive interpretation of the
free solvable group: this is simply the extension of the homology group of the Cayley graph of the
previous solvable group by the previous group $\Sol(n-1,d)$ itself with canonical cocycle. For the
metabelian group Sol(2,d) we obtain an especially simple interpretation: {\it the free metabelian
group $Sol(2,d)$ is the extension of the homology group of the lattice ${\Z}^d$ by the group $\Z^d$
with canonical 2-cocycle} - Theorem 4, (see \cite{Ver1}).

In the main part of the paper (paragraph 3) we describe the structure of metabelian group as
mentioned above extension of first homology group of Caley graph of lattice using canonical
cocycle. The simplest generators of the homology group is so called plaquettes (minimal cycles on
the lattice); for $d=2$ these are the free generators of the $H_1({\Z}^d)$, but for $d>2$ they are
not free and we describe the corresponding relations (lemma 5) \footnote{We are grateful to the
reviewer who pointed out that these relations are known  from Magnus imbedding}. In the theorem 4
and 5 we give the precise isomorphism between meatabelian groups those extensions.

The main technical point is the study of the sets of 2-cocycles of the group $\Z^d$ with values in
$H_1(\E^d)$ and the corresponding second cohomology group $H^2(\Z^d;H_1(\E^d))$. We prove that for
$d=2$ this group is isomorphic to $\Z$; its generator is the 2-cocycle that determines the
metabelian free group (Theorem 7); we also describe extensions that correspond to other 2-cocycles
(Theorems 8,9).

\bigskip

It seems that all these simple but fundamental facts were not mentioned in the literature and
text-books; at least the specialists in the theory of solvable groups and combinatorial group
theory whom the first author could asked about had not mentioned any references. But the goal of
all this investigation was to apply its to the studying of the representations the solvable groups
and asymptotic theory of random walks on such groups. In this paper we only mention one of these
questions - the normal form of the elements of the metabelian group which follows from homological
realization. Namely in the last paragraph we start from the geometrical description of the elements
of an arbitrary finitely generated group as a classes of pathes on the lattices. For metabelian
group this gives a convenient way to define the useful normal form which will be considered with
other applications of the constructed model, elsewhere (theorems 10,11).

We want to mention that the free solvable groups and their properties had been discussed in the
algebraic literature for a long time; see the books by Magnus et al.
 \cite{MagKarSol},Lyndon \cite{LySch}, Kurosh \cite{Kur}, articles by Malcev \cite{Malc}, Shmel'kin \cite{Shmel} on
the connection with wreath products, and so on. Some combinatorial and numerical information about
free solvable groups see in \cite{Sok, Eg, Petr}.

Remark that free metabelian groups are not a matrix groups because the commutant of it is
infinitely generated subgroup, - they are very close to the wreath products like wreath product of
infinite cyclic groups (see also \cite{Shmel}). This is why they are very useful in the theory of
random walks on the groups as it was considered on the examples of wreath products in the early
papers of the first author \cite{Ver3} and in the later papers of A.Ershler (see f.e.\cite{Er}).
Another important thing is that free solvable groups are in a sense {\it infinite dimensional
groups} - they could be considered as the lattices in some infinite-dimensional groups.
Consequently the representation theory of the free solvable groups (which is not well studied up to
now as we know) must be very instructive.

The main open problems concern to numerical and asymptotic characteristics of free solvable groups,
especially problems related to random walks on these groups, Poisson--Furstenberg boundaries
connections with wreath products, normal forms, and problems about representations of the
metabelian groups and its cohomology with values in the representations. It seems that the model
which we considered here will help to the solution of such questions.

\newpage
\section{Fundamental and homology groups of Cayley graphs}
\subsection{The fundamental group of the Cayley graph of a finitely generated group}

Let $G$ be a finitely generated group with a fixed set of generators $S=\{x_1, \ldots, x_d\}$. We
do not include the inverses of generators into the set $S$, but in what follows we will assume that
they are included into the alphabet.  Let $\Cal(G)$ be the right Cayley graph of the group $G$ in
these generators, i.e., the graph whose vertices are elements of $G$, with two vertices $g_1, g_2$
connected by a (directed) edge whenever $g_1^{-1}g_2 \in S$ is a generator\footnote{This means that
we multiply an element by some generator from the right.}

Let $F_d$ be the free group with $d$ free generators, which we will denote by $X_1, \ldots, X_d$.
There exists a normal subgroup $N\lhd F_d$ such that the sequence $$ 0\longrightarrow N
\stackrel{i}{\longrightarrow} F_d \stackrel{\phi}{\longrightarrow} G \longrightarrow 0
$$ is exact.

Here $i: N \rightarrow F_d$ is the inclusion map and the homomorphism $\phi:F_d\rightarrow G$ is
defined by the rule $X_i \rightarrow x_i$.

 Define Cayley graph $\Cal(G)$ as a one-dimensional topological complex, or as the space
 $G \times S \times [0,1]$ with the following identifications:

1) for a fixed $g$ and all $i=1,\dots, d$, all the points $(g,x_i,0)$ are identified with each
other and denoted by $(g,0)$;

2) the triples $(g,x_i,1)$ are identified with $(gx_i,0)$.

We denote this one-dimensional complex by the same symbol $\Cal(G)$.

The group $G$ acts on the graph $\Cal(G)$ from the left by shifts,  as a group isomorphism of the
Cayley graph. Namely, generator $x_i$ maps the vertex labelled by element $g\in G$ to the vertex
labelled by  $x_i g$. Note, that the free group $F_d$ also acts on the $\Cal(G)$: generator $X_i$
maps the vertex labelled by element $g\in G$ to the vertex labelled by  $x_i g$.

This actions induce a (left) action of $F_d$ on the fundamental group $\pi_1(\Cal(G))$ and a (left)
action of $G$ on the homology group $H_1(\Cal(G);\bf Z)$ of the complex $\Cal(G)$.

More precisely, let $\pi_1(\Cal(G))$ be the fundamental group with basepoint at the identity
element $e\in G$, i.e., the group of homotopy classes of loops on $\Cal(G)$ that start and end at
the point $e$.  An action of the group $F_d$ on $\pi_1(\Cal(G))$ is defined as follows. Let $X_i$
be one of the generators of $F_d$, and let $\bar \lambda\in \pi_1(\Cal(G))$ be the homotopy class
containing a loop $\lambda$. Then $X_i$ acts on $\bar \lambda$ according to the rule $$
X_i(\bar\lambda) = \bar\lambda^{X_i}, $$ where $\lambda^{X_i}$ is the loop that first goes from the
basepoint to the vertex $x_i\in \Cal(G)$ along the edge corresponding to $x_i$, then passes the
loop $\lambda$ shifted to this point, and then goes back to the basepoint along the same edge.

Now we will give a natural definition of a mapping $\gamma$ that
maps elements of the free group $F_d$ to directed paths on
$\Cal(G)$ starting at the basepoint.

Let $w\in F_d$ be some element of $F_d$. It has a unique representation as a
word in the alphabet $X_1, \ldots, X_d, X_1^{-1}, \ldots, X_d^{-1}$
that does not contain trivial cancellations.

\begin{defn} The path $\gamma (w)$ on $\Cal(G)$ is constructed
as follows: it starts at the point $e\in \Cal(G)$, goes
 with unit velocity along the edges corresponding to the
generators, and ends at the point $(w,0)\in \Cal(G)$.
\end{defn}

In a similar way, given a path that starts at the basepoint
and consists of edges of $\Cal(G)$, we can construct a word in the alphabet
$X_1, \ldots, X_d, X_1^{-1}, \ldots, X_d^{-1}$.

\begin{lem}
An element $w\in F_d$ belongs to the normal subgroup $N
\triangleleft F_d$ if and only if $\gamma(w)$ is a loop.
\end{lem}
\begin{proof}
Write an element $w$ as a word without cancellations:
$$
w=X_{k_1}^{\delta_1}X_{k_2}^{\delta_2}\cdots X_{k_n}^{\delta_n},
$$
where $\delta_i = \pm 1$.

This element belongs to $N$ if and only if
$x_{k_1}^{\delta_1}x_{k_2}^{\delta_2}\cdots x_{k_n}^{\delta_n}=e$
in $G$; this means that the last edge of $\gamma(a)$ connects
the vertex $x_{k_1}^{\delta_1}x_{k_2}^{\delta_2}\cdots
x_{k_{n-1}}^{\delta_{n-1}}$ with the basepoint, and thus  $\gamma(w)$
is a loop.
\end{proof}

Now we can define a mapping $\bar \gamma: N\rightarrow
\pi_1(\Cal(G))$ that maps an element $a\in N$ to the homotopy class
$\bar\gamma(a)$ of the loop $\gamma(a)$.

\begin{lem}
The mapping $\bar \gamma:N \rightarrow \pi_1(\Cal(G))$ is a group
isomorphism. Moreover, if $w\in F_d$ and $a\in N$, then
$$
\bar\gamma(a^w) =  w(\bar\gamma(a)),
$$
where $a^w = w a w^{-1}$.

\end{lem}
\begin{proof}
First of all, let us check that this mapping is a group
homomorphism. Let $u,v\in N$; then the path $\gamma (u)\gamma(v)$ is
the concatenation of the paths $\gamma (u)$ and $\gamma (v)$.
On the other hand, the word representing $uv$ differs from the
concatenation of the words representing $u$ and $v$ by trivial
cancellations, which lead to homotopically trivial subloops. Thus
the loop $\gamma(uv)$ is homotopic to $\gamma(u)\gamma(v)$, and
$\bar \gamma (uv)=\bar \gamma(u) \bar \gamma (v)$.

Note that every homotopy class from $\pi_1(\Cal(G))$ contains
at least one loop that consists of edges of $\Cal(G)$  and hence, by Lemma~1,
corresponds to some element of $N$. Thus $\bar \gamma$
is an epimorphism.

Now let us prove that if  $\gamma(a)$ is homotopic to the trivial loop,
then $a=e$. This follows immediately from the fact that $\Cal(G)$ is
an infinite bouquet of one-dimensional spheres and $N$ is a free
group, and, consequently, a homotopically trivial loop corresponds to the
identity element of $G$.

The condition $\bar\gamma(a^w) =  w(\bar\gamma(a))$, $a\in N$,
$w\in F_d$, follows directly from the definition of the action of
$F_d$ on $\pi_1(\Cal(G))$ and the construction of $\gamma$.
\end{proof}

Now we can deduce the following theorem.

\begin{thm}
The mapping $\bar \gamma$ defined above gives a canonical
$F_d$-equivariant isomorphism of the groups $N$ and
$\pi_1(\Cal(G))$. In other words, $\bar \gamma$ is an isomorphism
and
$$
\bar\gamma(a^w) =  w(\bar\gamma(a))
$$
for any $w\in F_d$ and $a\in N$.
\end{thm}

Note that both groups $N$ and $\pi_1(\Cal(G))$ are free, the first
one due to the Nilsen-Schreier theorem (see \cite{Kur}), and the
second one since $\Cal(G)$ is homotopically equivalent to a
bouquet of one-dimensional spheres. The theorem asserts that the
canonical isomorphism is equivivariant with respect to the action
of the free group.

\subsection{The homology group $H_1(\Cal(G);\bf Z)$}

As was mentioned in the previous section, the first homology group
$H_1(\Cal(G);\bf Z)$ can naturally be considered as a $G$-space.

Note that since $\Cal(G)$ is a one-dimensional complex, the group
$H_1(\Cal(G);\bf Z)$ is the group of one-dimensional cycles on
$\Cal(G)$, i.e., a subgroup of the free abelian group on all edges of
$\Cal(G)$ consisting of elements with trivial boundary, where the
boundary of an edge $((g,x_i,0),(g,x_i,1))$ is the null-dimensional
chain $(g,x_i,0)-(g,x_i,1)$.

The group $H_1(\Cal(G);\bf Z)$ is abelian, and we will write its
group law additively.

Note that $G$ acts on $H_1(\Cal(G);\bf Z)$
shifting  cycles by elements of $G$.

Let $N'$ be the commutant of the group $N$, i.e.,  the normal subgroup of
$N$ generated by commutators, elements of the form
$[x,y]=xyx^{-1}y^{-1}$, $x,y\in N$.

Note that if $N$ is a normal subgroup of $F_d$, then the
commutant $N'$ of N is also normal in $F^d$. Thus the group
$N/N'$ is a normal subgroup of $F_d/N'$, and the sequence
$$
0\longrightarrow N/N' \stackrel{i'}{\longrightarrow} F_d/N'
\stackrel{\phi'}{\longrightarrow} G \longrightarrow 0
$$
is exact. Here $i'$ is again the inclusion map and $\phi'$ is
defined by the rule $\pi(X_i)\mapsto x_i$, where $\pi:
F_d\rightarrow F_d/N'$ is the natural projection.

Note, that $N/N'$ has the natural structure of  a $G$-space
defined by the rule

$$
\pi(u)^{x_i} = \pi(X_i u X_i^{-1}),
$$
where $u\in N$. This action is well-defined. Indeed, suppose
$x_{i_1}^{\delta_1} x_{i_2}^{\delta_2} \ldots x_{i_k}^{\delta_k} =
e$. Then, $X_{i_1}^{\delta_1} X_{i_2}^{\delta_2} \ldots
X_{i_k}^{\delta_k} \in N$ and

$$
\pi(u)^{x_{i_1}^{\delta_1} x_{i_2}^{\delta_2}\ldots
x_{i_k}^{\delta_k}} = \pi(X_{i_1}^{\delta_1} X_{i_2}^{\delta_2}
\ldots X_{i_k}^{\delta_k})\pi(u)\pi(X_{i_1}^{\delta_1}
X_{i_2}^{\delta_2} \ldots X_{i_k}^{\delta_k})^{-1}=\pi(u),
$$
as $N/N'$ as abelian.

\begin{thm}
There is a $G$-equivariant isomorphism of the groups $N/N'$ and
$H_1(\Cal(G),\bf Z)$.
\end{thm}
\begin{proof}
Indeed, since $\Cal(G)$ is a linearly connected complex, the Gurevich theorem (see
\cite{FuchsRokh}) implies $$ \pi_1(\Cal(G))/(\pi_1(\Cal(G)))' \simeq H_1(\Cal(G);\bf Z). $$

Now the statement follows from the definition of the actions and from the Theorem 1.
\end{proof}

\subsection{The cocycle}

In what follows we will use the notation $$ U\leftthreetimes_{\tilde c} V $$ for the extension of a
group $V$ by a group $U$ determined by action of the group $U$ on the group $V$ as a group of
automorphisms, and by 2-cocycle $\tilde c$ of the group $U$ with values in the group $V$ which will
be described in each particular case.

The group $F_d/N'$ can be represented as an extension of the $G$-space
$H_1(\Cal(G);\bf Z)$ by the group $G$. Now we will describe
the cohomology class of cocycles corresponding to this extension.

For every element $g\in G$, fix a representation of $g$ as a word
in the alphabet of generators $\{x_1,\ldots, x_d\}=S$ and their
inverses. There are many such representations, and we will see that
different choices lead to cohomologous cocycles.

To fix a representation for every element means that for every $g\in
G$, we fix a path $\omega_g$ on $\Cal(G)$ connecting the
basepoint with the vertex $(g,0)$.

Let $g_1, g_2\in G$. Consider the loop $\lambda_{g_1,g_2}$ on
$\Cal(G)$ constructed as follows: $\lambda_{g_1,g_2}$
starts at the basepoint, then passes to the vertex $(g_1,0)$ along the path
$\omega_{g_1}$, then passes to the vertex $(g_1g_2,0)$ along the path
$g_1(\omega_{g_2})$ (i.e., the  path $\omega_{g_2}$ shifted by the element
$g_1$), and then returns to the basepoint along the path
$\omega_{g_1\cdot g_2}^{-1}$.

Denote by $\tilde \lambda_{g_1,g_2}$ the class in the homology
group $H_1(\Cal(G);\bf Z)$ that contains the loop
$\lambda_{g_1,g_2}$.

\begin{defn}Define a mapping $c:G\times G\rightarrow
H_1(\Cal(G);\bf Z)$ by the rule
$$
c(g_1,g_2)=\tilde \lambda_{g_1,g_2}.
$$
\end{defn}

\begin{lem} The mapping $c$ is a 2-cocycle of the group $G$ with values in $H_1(\Cal(G);\Z)$.
\end{lem}

\begin{proof}
We must prove that the mapping $c$ has trivial coboundary,
i.e., for any $g_1,g_2,g_3\in G$,
$$
(\partial c)(g_1,g_2,g_3) = c(g_1,g_2)+c(g_1\cdot g_2,g_3)-c(g_1,
g_2\cdot g_3)-g_1(c(g_2,g_3)) = 0.
$$
Indeed, the cycle $(\partial c)(g_1,g_2,g_3)$ corresponds to the
product of four loops $\lambda_{g_1,g_2}$, $\lambda_{g_1\cdot
g_2,g_3}$, $\lambda_{g_1,g_2\cdot g_3}^{-1}$, and
$g_1(\lambda_{g_2,g_3})^{-1}$, where
\begin{eqnarray*}
\lambda_{g_1,g_2} &=& \omega_{g_1} \cdot g_1(\omega_{g_2})\cdot
\omega_{g_1g_2}^{-1},
\\
\lambda_{g_1\cdot g_2,g_3}&=& \omega_{g_1g_2} \cdot
(g_1g_2)(\omega_{g_3})\cdot \omega_{g_1g_2g_3}^{-1},
\\
\lambda_{g_1,g_2\cdot g_3}^{-1} &=& \omega_{g_1g_2g_3} \cdot
g_1(\omega_{g_2g_3})^{-1}\cdot \omega_{g_1}^{-1},
\\
g_1(\lambda_{g_2,g_3})^{-1} &=& g_1(\omega_{g_2g_3})^{-1} \cdot g_1(
g_2(\omega_{g_3}))^{-1}\cdot g_1(\omega_{g_2} )^{-1}.
\end{eqnarray*}
(here the product of paths means their (coordinatewise) sum
as functions on the edges of $\Cal(G)$).

Note that $g_1(g_2(\omega_{g_3}))=(g_2g_1)(\omega_{g_3})$ and
each path occurs in the sum twice with opposite
directions. Hence $(\partial c)(g_1,g_2,g_3)=0$ in
$H_1(\Cal(G);\bf Z)$.
\end{proof}

Our definition of the cocycle $c$ depends on the choice of the set of
paths $\{\omega_g\}_{g\in G}$. The next lemma asserts that cocycles
constructed from different sets of paths are cohomologous.

\begin{lem}
Let $\{\omega_g\}_{g\in G}$ and $\{\omega'_g\}_{g\in G}$ be two
sets of paths. Then the corresponding cocycles $c= c(\{\omega_g\})$
and $c'=c'(\{\omega'_g\})$ are cohomologous.

Conversely, let $c= c(\{\omega_g\})$ be the cocycle constructed
from a set of paths $\{\omega_g\}_{g\in G}$, and let $c'$ be a
cocycle cohomologous to $c$. Then there exists a set of paths
$\{\omega_g\}_{g\in G}$ such that the cocycle $c' = c'(\{\omega_g\})$ is
constructed from $\{\omega_g\}_{g\in G}$.
\end{lem}
\begin{proof}
Cocycles $c$ and $c'$ are cohomologous if and only if there
exists a mapping $\alpha:G\rightarrow H_1(\Cal(G);\bf Z)$ such that
for any $g_1,g_2\in G$,
$$
c(g_1,g_2)-c'(g_1,g_2)= (\partial\alpha)(g_1,g_2)=
\alpha(g_1)+g_1(\alpha(g_2))-\alpha(g_1g_2).
$$

Now, if $c$ and $c'$ are constructed from the sets
$\{\omega_{g}\}_{g\in G}$ and $\{\omega'_{g}\}_{g\in G}$, then
the mapping $\alpha(g)$ that maps $g$ to the loop $\tilde \xi_g$, where
$\xi_g$ goes from the basepoint to $(g,0)$ along $\omega_g$ and then
returns to the basepoint along $(\omega'_g)^{-1}$ satisfies this condition.

To prove the converse claim, construct an appropriate set of
paths $\{\omega_g\}_{g\in G}$ using $\alpha$ and $\omega_g$.
\end{proof}

\begin{defn} Let $\tilde c$ be the cohomology class that consists
of all cocycles obtained from sets of paths
as described above.
\end{defn}

\begin{thm}
Let $F_d$ be the free group with $d$ generators, $N$ be some normal
subgroup of $F_d$, and $N'$ be the commutant of $N$. Let $\Cal(F_d/N)$
be the right Cayley graph of the quotient group $F_d/N$ associated with
the generators obtained
from the standard generators of $F_d$ through the canonical projection.

Then there is a canonical isomorphism
$$
F_d/N' \simeq F_d/N \leftthreetimes_{\tilde c} H_1(\Cal(F_d/N);\bf
Z),
$$
where $F_d/N$ acts on $H_1(\Cal(F_d/N);\bf Z)$ from the right by shifts and
the extension is determined by the cohomology class $\tilde c$
described above.
\end{thm}

\begin{proof}
Let $G=F_d/N$. It suffices to prove that the
cohomology class that determines $F_d/N'$ as an extension of
$N/N'$ by $G$ coincides with the class described above in this section.

Indeed, let $g\in G$. Denote by $u(g)\in F_d/N'$ some
representative of the class of $(F_d/N')/(N/N')$ corresponding to
$g$.

The cocycle $C$ that determines the extension is given by the formula
$$
u(g_1)u(g_2)=u(g_1g_2)C(g_1,g_2),
$$
where $C(g_1,g_2)\in N/N'$.

Fixing any representations of the elements $u(g)$ as words in
the alphabet $S\cup S^{-1}$ and considering the canonical image of
$C(g_1,g_2)$ in $H_1(\Cal(G);\bf Z)$, we arrive at the cocycle $c$
constructed from the set of paths determined by the elements
$\{u(g)\}_{g\in G}$.
\end{proof}

In the next section we will consider the special case when
$N=F_d^{(n-1)}$ and $N' = F_d^{(n)}$.

\newpage

\section{Free solvable groups}

\subsection{The free solvable group $\Sol(n,d)$}

In this section we will apply Theorem 2 to the special case
when $N$ is the $(n-1)$th commutant of the free group $F_d$.

\begin{defn}
Let $\Sol(n,d)$ be the free solvable group of degree $n$ with
$d$ generators, i.e., the free object in the category of all solvable
groups of degree $n$ with $d$ generators.

The group $\Sol(n,d)$ is isomorphic to the quotient of the free
group with $d$ generators $F_d$ over its $n$th commutant
$F_d^{(n)}$:
$$
\Sol(n,d)\simeq F_d/F_d^{(n)}.
$$

Let $\Met(d)=\Sol(2,d)$ be the free metabelian group, i.e., the
free solvable group of degree 2 with $d$ generators.
\end{defn}

Let us take the symmetric set of generators $S=\{x_1, \ldots, x_d
\}$ in $\Sol(n,d)$, where $x_i$ corresponds to the standard generator
 $X_i$ of the free group.

Applying Theorem 3 to $N=F_d^{(n-1)}$, we arrive at the following
theorem.

\begin{thm}
There is a canonical isomorphism
$$
\Sol(n,d) \simeq \Sol(n-1,d) \leftthreetimes_{\tilde c}
H_1(\Cal(\Sol(n-1,d))),
$$
where the extension in the right-hand side is determined by the right action of
$\Sol(n-1,d)$ on $H_1(\Cal(\Sol(n-1,d)))$ by shifts and by the
cohomology class $\tilde c$ defined in the same way as in
Section 2.3.
\end{thm}
\begin{proof}
Indeed, this theorem is an immediate consequence of Theorem 3 with
$N=F_d^{(n-1)}$. In this case $N'=F_d^{(n)}$ and $G \simeq
F_d/F_d^{(n)}\simeq \Sol(n-1,d)$.
\end{proof}

In particular, consider the case $n=2$.

Note that $\Sol(1,d)$ is the free abelian group ${\bf Z}^d$ and
its Cayley graph $\Cal({\bf Z}^d)$ is the $d$-dimensional integer
lattice in ${\bf R}^d$, i.e., the vertices of $\Cal({\bf Z}^d)$ are
points of ${\bf R}^d$ with integer coordinates and its edges are unit
edges connecting these points. As a  one-dimensional complex,
$\Cal({\bf Z}^d)$ is the union of all lines in ${\bf R}^d$ parallel
to axes and crossing integer points. Let us denote this complex by
${\bf E}^d$.

As an immediate consequence of the previous theorem, we obtain the following
theorem, which was stated in \cite{Ver1}.

\begin{thm} There is a canonical isomorphism
$$
\Met(d) \simeq {\bf Z}^d \leftthreetimes_{\tilde c} H_1({\bf
E}^d),
$$
where $\Z^d$ acts on cycles by shifts from the right and the cohomology
class is defined as before.
\end{thm}

In the next section we will describe this isomorphism more precisely.

\subsection{The free metabelian group $\Met(d)$}

Let $\bar m= (m_1, \ldots, m_d)\in {\bf Z}^d$ be an integer vector
and $ 1\leq i,j\leq d$. Consider the unit square $P(\bar m, i, j)$ in
${\bf R}^d$ with coordinates of vertices
$$
\begin{array}{c} (m_1,\ldots,m_i,\ldots,m_j,\ldots,m_d), \\
(m_1,\ldots,m_i+1,\ldots,m_j,\ldots,m_d), \\
(m_1,\ldots,m_i+1,\ldots,m_j+1,\ldots,m_d),\\
(m_1,\ldots,m_i,\ldots, m_j+1,\ldots,m_d),
\end{array}
$$ i.e., $P(\bar m,i,j)$ lies in the plane parallel to the
coordinate plane defined by the $i$th and $j$th axes, and the
point of $P(\bar m,i,j)$
closest to the origin  has coordinates
$\bar m$.

\begin{defn}Let $p(\bar m, i,j)\in H_1(\E^d))$ be the oriented boundary
of $P(\bar m,i,j)$, where the orientation is induced by the order of
coordinates $i,j$. Elements of these form will be called {\it
plaquettes}.
\end{defn}

\begin{lem}
The group $H_1(\E^d)$ is generated by plaquettes $p(\bar m,i,j)$,
where $\bar m\in {\bf Z}$, $1\leq i<j\leq d$.

For $d=2$, plaquettes generate the group $H_1(\E^d)$ freely.

For $d\geq 3$, the set of relations on plaquettes is generated by
relations of the form
$$
\begin{array}{c}
p(\bar m,i,j)+p(\bar m,i,k)+p(\bar m,j,k)\\
-p(\bar m+e_k,i,j)-p(\bar m+e_j,i,k)+p(\bar m+e_i,j,k)=0,
\end{array}
$$
where $1\leq i <j<k\leq d$, $\bar m\in \Z^d$, and $e_i, e_j, e_k$
are unit coordinate vectors of the corresponding axes.

In other words, the relations on the set of generators
$\{p(\bar m ,i, j)\}$
are generated by the ``cube'' relation: for every unit
3-dimensional cube with integer vertices, the sum of the oriented
boundaries of its six faces is equal to zero.
\end{lem}

\begin{proof}
Consider the cell structure on ${\bf R}^d$ with $k$-dimensional
cells being unit $k$-dimensional cubes with integer vertices. Let
$C_k$ be the free abelian group on the set of all $k$-dimensional
cells. Defining the boundary $\partial$ in a standard way, we have the
exact sequence
$$
0\stackrel{\partial}{\longleftarrow}
C_0\stackrel{\partial}{\longleftarrow}
C_1\stackrel{\partial}{\longleftarrow} \cdots
\stackrel{\partial}{\longleftarrow}
C_d\stackrel{\partial}{\longleftarrow} 0.
$$
Note that $H_1(\E^d)\simeq Z_1(\E^d),$ where $Z_1(\E^d)$ is the
kernel of the boundary in $C_1$. 

Since this sequence is exact,
$$
H_1(\E^d)\simeq C_2/\partial(C_3),
$$
but $C_2$ is generated by the elements $P(\bar m ,i,j)$, $\bar m \in
\Z^d$, $1\leq i,j\leq d$. Hence $H_1(\E^d)$ is generated by the
boundaries of these elements, which are plaquettes. The relations are
defined by elements of $\partial(C_3)$ that are generated by the sums
of faces of unit 3-dimensional cubes. In the case $d=2$, we have
$\partial(C_3)=0$, and hence plaquettes are free generators.

The theorem can also be proved directly using induction and simple
geometrical considerations.
\end{proof}

Now let $\p(\bar m ,i,j)\in \Met(d)$ be the element of the free
metabelian group defined by the formula
$$
\p(\bar m, i, j)= x_1^{m_1}x_2^{m_2}\cdots x_d^{m_d} \cdot
[x_i,x_j]x_d^{-m_d}\cdots x_2^{-m_2}x_1^{-m_1}.
$$

Fix the set $\{\omega_{\bar m}\}_{\bar m\in \Z^d}$ of paths
defined by the formula
$$
\omega_{\bar m}=\gamma(X_1^{m_1}X_2^{m_2}\cdots X_d^{m_d})
$$
(i.e., the path $\omega_{\bar m}$ first goes along the axis $X_1$,
then along $X_2$, etc., until it reaches the point $\bar m$).

Consider the cocycle $c:\Z^d\times \Z^d \rightarrow H_1({\bf
E}^d)$ defined as in Section~1.3 with this choice of
$\omega_{\bar m}$. Let $\tilde c$ be the cohomology class
containing $c$.

Consider the extension ${\bf Z}^d \leftthreetimes_{\tilde c}
H_1({\bf E}^d)$ determined by the
action of $\Z^d$ on plaquettes defined by the formula
$$
\bar m(p(\bar n, i,j))=p(\bar n+ \bar m, i,j)
$$
and the class $\tilde c$. Obviously, this extension is generated by
elements of the form $(e_k,
p(\bar m,i,j))$, $\bar m\in \Z^d$, $k, i<j\in Z$.

\begin{thm}
There is a canonical isomorphism
$$\Met(d) \simeq {\bf Z}^d \leftthreetimes_{\tilde c} H_1({\bf
E}^d) $$
defined on the generators by the rule
$$
\alpha((e_k, p(\bar m,i,j))) = x_k \p(\bar m,i,j).
$$
\end{thm}

\begin{proof}
The theorem follows directly from Theorem 4.
\end{proof}

\subsection{The cohomology group $H^2(\Z^2, H_1(\E^2))$}

Let $H^2(\Z^2, H_1(\E^2))$ be the cohomology group with respect to
a fixed right action of $\Z^2$ on $H_1(\E^2)$ by shifts.

In this section we will prove that this group  is isomorphic to
the group of integers $\Z$, and under this isomorphism, the
class $\tilde c$ that determines the free metabelian group $\Met(2)$
corresponds to the identity element of $\Z$.

Let $\tilde y\in H^2(\Z^2, H_1(\E^2))$ be a cohomology class,
and let $y:\Z^2\times\Z^2\rightarrow H_1(\E^2)$ be some
representative of this class.

Let $G$ be the extension of $H_1(\E^2)$ by $\Z^2$ determined by
the class $\tilde y$; in particular, we consider the
representation determined by the cocycle $y$. In other words,
every element of $G$ can be written as a pair $((m,n), \gamma)$,
where $m,n\in \Z$ and $\gamma\in H_1(\E^2)$ is a cycle.

The group law in $G$ is determined by the formula
$$
((m,n),\gamma_1)\cdot((k,l),\gamma_2) = ((m+k,n+l),
\gamma_1^{(k,l)}+\gamma_2+y((m,n),(k,l))),
$$
where, as before, $\gamma_1^{(k,l)}$ denotes the cycle $\gamma_1$
shifted by the vector $(k,l)$.

Denote the vector $(0,1)$ by $e_1$ and the vector $(1,0)$ by $e_2$.
We will write $0$
for the origin $(0,0)\in \Z^2$.

In every cohomology class there is at least one normalized
element; thus we may assume that the cocycle $y$ is normalized, i.e.,
$y(0, (m,n)) = y((m,n),0)= 0$ for every integer vector $(m,n)$.

Consider the commutator of the elements $(e_1,0)$ and $(e_2,0)$ in the
group $G$:
\begin{eqnarray*}
[(e_1,0),(e_2,0)] &=&
(e_1,0)(e_2,0)(-e_1,-y(e_1,-e_1))(-e_2,-y(e_2,-e_2))
\\
&=&(0, y(e_1,e_2)^{-e_1-e_2} -
y(e_1,-e_1)^{-e_2}+y(e_1+e_2,-e_1)^{-e_2} ).
\end{eqnarray*}

Let
$$
\gamma_y = y(e_1,e_2)^{-e_1-e_2} -
y(e_1,-e_1)^{-e_2}+y(e_1+e_2,-e_1)^{-e_2}.
$$

\begin{defn}Set
$$
\beta(y) = S(\gamma_y),
$$
where $S(\gamma_y)$ is the algebraic area of $\gamma_y$, i.e., the
sum of the winding numbers  of all integer
points in the domain bounded by $\gamma_y$ with respect to $\gamma_y$.
\end{defn}

As follows from Lemma 5, the cycle $\gamma_y$ has a unique
representation as a sum of plaquettes:
$$ \gamma_y = \sum\limits_{(m,n)} k_{(m,n)}p(m,n), $$
where the sum is
taken over all integer points $(m,n)\in \Z^2$,  $k_{(m,n)}\in Z$
are integer coefficients (only finitely many of them are not
equal to zero), and $p(m,n)=p((m,n),1,2)$ are plaquettes as defined in
Section~3.2.

Then
$$\beta(y) = S(\gamma_y)= \sum\limits_{(m,n)}k_{(m,n)}.$$

The latter equality can also be taken as a definition of
$S(\gamma_y)$.

\begin{thm} The constructed mapping $\beta$ is constant on
cohomology classes, i.e., if cocycles $y$ and $y'$ are
cohomologous, then $\beta(y)=\beta(y')$.

The well-defined homomorphism
$$
\beta:H^2(\Z^2, H_1(\E^2))\simeq \Z
$$
is a group isomorphism, and
$$ \beta(\tilde c)=1, $$
where the class $\tilde c$ determines the
extension isomorphic to the free metabelian group $\Met(2)$.
\end{thm}
\begin{proof}
Let $y'\in \tilde y$. We have to show that
$S(\gamma_{y'})=S(\gamma_y)$.

Indeed, since $y$ and $y_1$ are cohomologous, there exists a mapping
$u:\Z^2\rightarrow H_1(\E^2)$ such that for any $(m,n)$ and
$(k,l)$ from $\Z^2$,
$$
y'((m,n),(k,l)) = y((m,n),(k,l)) + u(k,l)+u(m,n)^{(k,l)}-
u(m+k,n+l).
$$

Then
$$ \gamma_{y'} = \gamma_y + u(e_1)^{e_2-e_1} - u(e_1)^{-e_1} +
u(e_2)^{-e_1} - u(e_2).\eqno{(*)} $$

It is obvious that the algebraic area of the difference
$\gamma_{y'}-\gamma_y$ is equal to zero.

By definition, $\beta$ is a group homomorphism. Now we will show
that $\beta$ is an epimorphism. Note that if $\tilde c $ is the cocycle
defined in Section~2.1, then $\gamma_c = p(0,-1)$, and thus
$\beta(\tilde c) = 1$. Then for every $n\in \Z$ we have
$$
\beta(\tilde{nc}) = n,
$$
and $\beta$ is an epimorphism.

It suffices to show that $\beta$ is a monomorphism. Assume that
$\tilde y\in H_1(\E^2)$ is a cohomology class, $y\in \tilde y$ is its
representative, and $\beta(\tilde y)=S(\gamma_y)= 0$.

If $S(\gamma_y)=0$, then, using equation (*), we can find a mapping
$u:\Z^2\rightarrow H_1(\E^2)$ such that $\gamma_{y'}=0$. This means
that in the extension $G$ of $H_1(\E^2)$ by $\Z^2$
determined by $y'$, the elements $(e_1,0)$
and $(e_2,0)$ commute.

Let $H\subset G$ be the subgroup of $G$ generated by the elements
$(e_1,0),(e_2,0)$, and let $N\lhd G$ be the normal subgroup of $G$
consisting of all elements of the form $(0,\gamma)$, $\gamma\in
H_1(\E^2)$. Then $H\simeq \Z^2$, $N\simeq H_1(\E^2)$.

Obviously, $H\cap N = \emptyset$ and $HN = G$. Then
$G$ is a split extension,
and the mapping $y$ is zero.
\end{proof}

\subsection{The structure of other extensions}

 As was proved in the previous section, the cohomology
group $H^2(\Z^2, H_1(\E^2))$ is isomorphic to the group $\Z$, and
the free metabelian group $\Met(2)$ corresponds to the
generator $1 \in \Z$. It is
interesting to see what is the structure of extensions
corresponding to other integers.

Let $k\in\Z$ be an integer. Consider the group
$\Met_k(2)$ constructed as follows.

The group $\Met_k$(2) has three generators $x,y,z$. Let
$$
V_g h =  g h g^{-1}, \qquad g,h\in \Met_k(2).
$$
Then the group $\Met_k(2)$ is determined by the relations
\begin{eqnarray}
[x,y]&=&z^k,\\
\ \ \  [V_{x^m y^n}z, V_{x^{m'} y^{n'}}z] &=& [V_{x^{m'} y^{n'}}z,
V_{x^m y^n}z]\qquad\forall m,n\in \Z.
\end{eqnarray}

Note that $\Met_1(2)=\Met(2)$.

Consider the normal subgroup $N$ of $\Met_k(2)$ generated by
the element $z$ (in the sense that $N$ is the smallest normal subgroup of
$\Met_k(2)$ that contains $z$); in other words, $N$ consists of
elements of the form $\{z^{x^my^n}\}_{m,n\in \Z}$ and their
products.

Note that for $k=\pm 1$, the subgroup $N$ is equal to the
commutant $\Met_k(2)'$, but for other values of $k$, the commutant
is strictly included into $N$.

\begin{thm}
1) There is a group isomorphism $\Met_k(2)/N\simeq \Z^2$.

2) The group $N$ is isomorphic to $H_1(\E^2;\Z)$ as a group with
operators from $\Z^2$. 

3) Let $\tilde c_k\in H^2(\Z^2,H_1(\E^2))$ be the class of
2-cocycles that determines the group $\Met_k(2)$ as an extension
of the group $N\simeq H_1(\E^2)$ by the group $\Met_k(2)/N\simeq
\Z^2$. Then
$$ \beta(\tilde c_k)=k, $$
where the mapping $\beta$ is defined in Section~3.3.

\end{thm}
\begin{proof}
1) This follows from the fact that $[x,y]=e \mod N$ and there are
no other relations on the elements $x$ and $y$ modulo $N$.

2) Indeed, the isomorphism is given by the correspondence
$$ z^{x^my^n}\mapsto p(m,n), $$
where $p(m,n)\in H_1(\E^2)$ is the plaquette defined in the previous
section.

3) Let $c_k\in \tilde c_k$ be the cocycle corresponding to the
following set of representatives:
$$ (m,n)\mapsto x^m y^n. $$
Then the value of the cocycle $c_k((m,n),(m',n'))$ is a cycle 
corresponding to the word $x^my^nx^{m'}y^{n'}y^{-n-n'}x^{-m-m'}$.

Thus $c_k((1,0),(1,0))=0$, $c_k((1,0),(0,1))=0$,  and the value of 
$c_k((1,1),(-1,0))$ corresponds to the word $[x,y]=z^k$, i.e,
$$ c_k((1,1),(-1,0))=k p(0,0), $$
and
\begin{eqnarray*}
\gamma_{c_k} &=& c_k((1,0),(0,1))^{(-1,-1)} -
c_k((1,0),(-1,0))^{(0,-1)}+c_k((1,1),(-1,0))^{(0,-1)}\\
&=&kp(0,-1);
\end{eqnarray*}
thus
$$ \beta(\tilde c_k) = S(\gamma_{c_k}) = k. $$
\end{proof}
Remeber that set of placuettes $p(n,m)$ is the set of free generators
of $H_1(E^2)$, so each cycle is a finite sum $h=\sum k(n,m)p(n,m)$.
Let $M$ be the normal subgroup of $N$ generated by $z^k$. In other
words, $M$ consists of elements of the form $z^{k x^m y^n}$,
$m,n\in\Z$, and their products. Theorem 8 implies the following
assertions.

\begin{thm}
1) The groups $N$ and $M$ are isomorphic to $H_1(\E^2)$. Let us
fix isomorphism between $N$ and $H_1(\E^2)$ then the group $M$
goes to subgroup of $H_1(\E^2)$ generated by all the cycles
multiplicities of which with respect to each plaquettes belong to
the group $k\Z$,

2) 2-cocycles of the class $\tilde c_k$ have values in the group
$M$ as subgroup of $H_1(\E^2)$  (see item 1))

3)$$ N/M \simeq \oplus_{\Z^2} \Z_{|k|},$$
$$ \Met_k(2)/M \simeq \Z^2 \leftthreetimes
\oplus_{\Z^2}\Z_{|k|},$$
where $\Z_{|k|}$ is the finite cyclic group of order $|k|$.
and $\leftthreetimes$ denotes the semidirect product with
respect to the action of $\Z^2$ on $\oplus_{\Z^2}\Z_{|k|}$ by
shifts of coordinates.

4) The commutant $\Met_k(2)'$ of the group
$\Met_k(2)$ is subgroup of $N\subset H_1(\E^2)$ and has the following
description:
$$M\subset \Met_k(2)'=\{h=\sum k(n,m)p(n,m) \in H_1(\E^2): \sum_{n,m}k(n,m) \in k\Z\}\subset N$$

so the corresponding quotient is the following:
 $$ \Met_k(2)/\Met_k(2)'=\Z^2\oplus \Z_{|k|}.$$
\end{thm}
\begin{proof}
The first three statements are obvious. The last one follows from
the fact, that

$$
[x,z] = V_x z \cdot z^{-1} = p(0,1)-p(0,0),$$$$ [y,z] = V_y z
\cdot z^{-1} = p(1,0) - p(0,0),$$$$ [x,y] = z^k = k p(0,0)
$$

and considering other possible combinations of commutators we can
get any element with algebraical area being in $k\Z$.
\end{proof}

\newpage
\section{Geometrical description of a group as a set of classes of
paths on the lattice ${\Z}^d$}

We give a useful geometrical model of the elements of an arbitrary
countable finitely generated group $G$ as a class $_G$-equivalent
pathes on the lattice ${\Z}^d$ (see \cite{Ver1}).

 Let $F_d$ be a free group with standard generators $X_1 \ldots X_d$.
 Consider the free abelian group $\Z^d$. Let $\E^d$ be a Cayley graph of the
group $\Z^d$, in other words, $\E^d$ consists of all unit edges,
connecting integer points in $\R^d$.

Let

$$
w = X_{i_1}^{\delta_1} X_{i_2}^{\delta_2}\ldots X_{i_k}^{\delta_k}
\in F_d
$$
be some word. We can map this word to the directed path on the
lattice $\E^d$ starting at the origin the same way as we did it in
paragraph 2.1. Namely, label each of the positive orts with
integers $1,\ldots,d$ and start reading the word $w$ from left to
right. Meeting the letter $X_{i_l}^{\delta_l}$ we add to the path
one unit edge, starting at the end point of the previous one and
going in the direction of the $i_l$th ort if $\delta_l=1$ and in
the opposite direction otherwise.

Lets consider the set $\Pi_d$ of all directed paths on the lattice
$\E^d$ starting at the origin and consisting of unit edges. If two
paths differ only by edge which is included into one of the paths
in one direction and on the next step in the opposite one, we will
consider these paths to be the same. Thus, the path corresponding
to the word of the form $ww^{-1}$ is considered to be empty and we
so have one-to-one correspondence between elements of the group
$F_d$ and paths from the set $\Pi_d$.

Note, that concatenation of paths makes $\Pi_d$ into the group,
isomorphic to $F_d$.

Now, consider some finitely generated group $G$ with a fixed set
of generators $S=\{x_1, \ldots ,x_d\}$. This group can be
presented as a quotient group of the free group $F_d$:

$$
G= F_d/N,
$$
where $N\lhd F_d$ and the canonical projection $\pi:F_d\rightarrow
G$ maps $X_i$ to $x_i$.

Every element $g$ of the group $G$ corresponds to a class of words
of the group $F_d$. This class consists of all possible words in
the alphabet of $X_i, i=1,\ldots,d$ and their reverses such, that,
these words with $X_i$ replaced by $x_i$ represent element $g$.

As every element of the free group corresponds to the unique path
from $\Pi_d$, elements of the group $G$ correspond to classes of
paths of $\Pi_d$. Let's say two paths to be $G$-equivalent iff
they belong to the same class, corresponding to some element of
$G$.

For some groups it is possible to give independent geometrical
description of this equivalence relation.

\noindent\bf Example 1.\rm The free abelian group $\Z^d$. In this
case the class of paths corresponding to element $\bar m\in \Z^d$
is obviously defined by the end point: all the paths from this
class end at the point $\bar m$. In other words, two paths are
$\Z^d$-equivalent iff they have the same end-point.

\noindent\bf Example 2.\rm Free two-step nilpotent groups
(=discrete Heizenberg group). It is a less trivial example. Let
$G$ be the free two-step nilpotent group with $d$ generators.
First we consider the case $d=2$; this is the discrete Heisenberg
group, that is the group of integer upper triangular matrixes with
1s on the main diagonal. We denotes the generators by $a$ and $b$.
Then the relations are as follows:

$$
[a,b]a = a[a,b], \ \ [a,b]b = b[a,b], \ \ [a,b] = a b
a^{-1}b^{-1}.
$$
\begin{thm} In  the case of discrete Heisenberg group $G=Heis$
a closed pathes is $G$-equivalent to the zero path if and only if
algebraic (oriented) area enclosed by this path is zero. The two
arbitrary paths $\gamma_1$ and $\gamma_2$ are $G$-equivalent, iff
the (closed) path $\gamma_1 \gamma_2^{-1}$ is equivalent to the
zero path.
\end{thm}
For $d>2$ (free two-step nilpotent group with $d$ generators) the
equivalence of a path to the zero path is described by the same
criterion applied for all projections of the closed path on the
2-dimensional coordinate subspaces of the lattice

Now following to \cite{Ver1} we will describe the equivalence of
the paths on the lattice ${\Z}^d$ for the free metabelian group
$\Met(d)$.

First of all define the {\it algebraic multiplicity of the edge in
the given path}. Consider the free metabelian group $\Met(d)$ with
$d$ generators. Let $\gamma$ be some path from $\Pi_d$ and let
$e\in \E^d$ be some edge.

Then $n_{\gamma}(e)$ to be the difference between number of times
$e$ is passed by $\gamma$ in positive direction and number of
times $e$ is passed by $\gamma$.

\begin{thm} Pathes $\gamma_1,\gamma_2\in \Pi_d$ are
$Met(d)$-equivalent if and only if for all $e\in \E^d$

$$
n_{\gamma_1}(e) = n_{\gamma_2}(e).
$$

\end{thm}
\begin{proof}
This theorem directely follow from the theorem 6.
\end{proof}

\end{document}